\documentclass[a4paper]{amsart}
\usepackage{graphicx}

\newtheorem{thm}{Theorem}[section]
\newtheorem{cor}[thm]{Corollary}
\newtheorem{lema}[thm]{Lemma}
\newtheorem{prop}[thm]{Proposition}
\newtheorem{exam}[thm]{Example}

\theoremstyle{definition}
\newtheorem{defn}[thm]{Definition}

\theoremstyle{remark}
\newtheorem{rem}[thm]{Remark}

\numberwithin{equation}{section}

\newcommand{\R}{\mathbb R}
\newcommand{\N}{\mathbb N}

\newcommand{\ve}{\varepsilon}

\newcommand{\lam}{\lambda}

\begin{document}
\title[Refined asymptotics]{Refined asymptotics for eigenvalues on domains of infinite measure}

\author[J. Fern\'andez Bonder, J.P. Pinasco and A.M. Salort]{Juli\'an Fern\'andez Bonder, Juan Pablo Pinasco and Ariel M. Salort}

\address{Departamento de Matem\'atica, FCEyN - Universidad de Buenos Aires\hfill\break \indent Ciudad Universitaria, Pabell\'on I  (1428) Buenos Aires,
Argentina.}

\email[J. Fernandez Bonder]{jfbonder@dm.uba.ar}

\urladdr[J. Fernandez Bonder]{http://mate.dm.uba.ar/~jfbonder}

\email[J.P. Pinasco]{jpinasco@dm.uba.ar}

\urladdr[J.P. Pinasco]{http://mate.dm.uba.ar/~jpinasco}

\email[A.M. Salort]{asalort@dm.uba.ar}



\begin{abstract}
In this work we study the asymptotic distribution of eigenvalues in one-dimensional open sets. The method of proof is rather elementary, based on the Dirichlet lattice points problem, which enable us to consider sets with infinite measure. Also, we derive some estimates for the the spectral counting function of the Laplace operator on unbounded two-dimensional domains.
\end{abstract}

\maketitle

\section{Introduction}

Let us consider an open set $\Omega \subset \R$ which is a disjoint union of bounded intervals, $\Omega=\bigcup_{j\in\N}I_j$. Let us suppose that the lengths of the intervals are decreasing and goes to zero,
$$
|I_1| \ge |I_2| \ge \cdots \ge |I_j|\ge \cdots \searrow 0.
$$
We can assume that there exists some nonincreasing function $g:(0, \infty) \to (0,\infty)$ such that
$$
|I_j| = g(j).
$$
Now, we may consider the following problems:

\begin{itemize}
\item A Lattice Point Problem: to estimate, for $x\nearrow \infty$, the number of lattice points below the curve $x g(t)$,
\begin{equation} \label{contarpuntos} 
N(x) = \# \{(j, k) \in \N \times \N \colon k \le x g(j) \} = \sum_{j=1}^\infty [ x g(j)].
\end{equation}

\item An Eigenvalue Counting Problem: to estimate, for $\lam \nearrow \infty$, the number of eigenvalues less than or equal to $\lam$ of $-u''=\lam u$ in $\Omega$
with zero Dirichlet boundary conditions on $\partial \Omega$,
$$
N(\lam) = \# \{j\in\N\colon \lam_j \le \lam\},
$$
\end{itemize}

The first one is called a {\it plane multiplicative problem}, following Kr\"atzel \cite{KRA}, and generalizes the {\em Dirichlet's divisor problem}, that is, to count the asymptotic number of divisors of the integers less than or equal to $x$, which is equivalent to count the number of lattice points below the hyperbola $y = x/t$ in the first quadrant.

The second one is a one dimensional variant of an old problem, {\em Can one hear the dimension of a drum?}, to use the catchy description given to it by Kac \cite{Kac}. The idea behind this name is the following: the square root of the eigenvalues of the Laplace operator in $\Omega \subset \R^2$ coincide with the musical notes of a membrane with the shape of $\Omega$, and we can ask about the geometric properties of $\Omega$ which can be inferred from the sequence of eigenvalues Here, we are interested in the dimension of the boundary of a {\em fractal string} $\Omega$, as Lapidus called this kind of sets \cite{La1}.

Indeed, both problems are the same: the eigenvalues of $-u'' = \lam u$ in $I_j$ are $\{\frac{\pi^2 k^2}{g(j)^2}\}_{k\ge 1}$, and we have
\begin{equation}\label{autovalores}
\begin{aligned}
N(\lam)  &=\sum_{j=1}^{\infty} \#\Big\{k\in\N\colon \frac{\pi^2 k^2}{g(j)^2} \le \lam \Big\}\\
&= \sum_{j=1}^{\infty} \#\Big\{k\in\N\colon k \le\frac{g(j)\lam^{1/2}}{\pi} \Big\} \\
&= \sum_{j=1}^{\infty} \Big[ \frac{g(j)\lam^{1/2}}{\pi} \Big]
\end{aligned}
\end{equation}
So, calling $x = \frac{\lam^{1/2}}{\pi}$, this expression coincides with equation \eqref{contarpuntos}, and we see that there exists a connection between the Dirichlet problem and the asymptotic behavior of eigenvalues. Let us mention that the eigenvalue counting problem for the Laplacian when $\Omega$ is the unit square in $\R^2$ coincide with the Gauss Circle Problem, i.e., to estimate the number of lattice points inside an expanding circle (see \cite{Duke}).

In this work we are interested in the asymptotic number of eigenvalues of the following $p-$laplacian eigenvalue problem in $\Omega$:
\begin{equation} \label{ecu}
-(|u'|^{p-2}u')'=\lam |u|^{p-2}u,
\end{equation}
with zero Dirichlet boundary conditions on $\partial \Omega$, and $1<p<+\infty$. Beside some technical details (see Section \S 2), there exists a closed formula for the $p-$laplacian eigenvalues similar to the one of the linear problem, which gives the full spectrum (see \cite{MANA, BOPI}). In this sense, our paper is only a minor generalization of previous works of He, Lapidus and Pomerance \cite{HL, LA} where the case $p=2$ was considered, and whenever the measure of $\Omega$ is finite, we obtain that
$$
N(\lam) = \#\{j\in\N\colon \lam_j \le \lam \} = \frac{|\Omega|}{\pi_p} \lam^{1/p} + \frac{\zeta(d)}{\pi_p^d} f(\lam^{1/p}) + o(f(\lam^{1/p}))
$$
by replacing everywhere $\lam^{1/2}$ by $\lam^{1/p}$, and $\pi$ by $\pi_p$ in the results of \cite{HL}. In the previous formula, $f(\lam^{1/p}) = g^{-1}(\lam^{-1/p})$, for $0<d<1$, and $\zeta$ is the Riemann Zeta function. The term $f(\lam^{1/p})$ is connected with a generalized notion of fractal dimension, and we have $f(\lam^{1/p}) = \lam^{d/p}$ when the Minkowski dimension of $\partial \Omega$ is $d$. The precise definitions and properties of $g$ and related functions is given in Section \S 3, together with the definitions of the generalized Minkowski content and dimension.

The proofs in those works depends on difficult estimates of the remainder terms of certain convergent series. We present in Section \S 4 a simplified proof based on the equivalence of the two problems stated above and some arguments from number theory. When the lengths of the intervals satisfy $|I_j| \sim j^{-1/d}$, as in \cite{LA}, this ideas were used in \cite{PINA}.

However, as a by-product of the number theoretic methods, we are able to extend those results to fractal strings $\Omega$ with infinite measure, and this is the main aim of our work. Let us observe that the sum in equation \eqref{autovalores} is well defined whenever $g(t) \searrow 0$ for $t\nearrow\infty$, even when $\sum_{j=1}^{\infty} g(j)$ diverges.

So, in Section \S 5 we characterize the growth of the number of eigenvalues $N(\lam)$ in terms of the decay of the lengths of the intervals when the measure of $\Omega$ is not finite. We obtain the following non-standard asymptotic formula
$$
N(\lam) = \#\{j\in\N\colon \lam_j \le \lam \} = \frac{\zeta(d)}{\pi_p^d} f(\lam^{1/p}) + o(f(\lam^{1/p})),
$$
where now $d>1$.

In the finite measure case, the term depending on $f$ can be thought as a boundary contribution. The measure of $\Omega$ gives the main term of the asymptotic of the number of lattice points, and the second term can be understood as the number of points which are close to the boundary and enter when we dilate slightly the domain. Now, when the measure of $\Omega$ is infinite, the main term is still a boundary term, which shows the asymptotic growth of the measure of the domain; in this case, when we dilate slightly the domain, a huge number of lattice points enter although it has exactly the same form that the second term in the other case.

The discreteness of the spectrum of an elliptic operator is not well understood yet when the domain has infinite measure. We refer the interested reader to \cite{BL, CH, Hew, Hw, BS} where a special class of sets in $\R^N$ is considered ({\it horn-shaped domains}, a $N-1$ dimensional set scaled in the other dimension). In \cite{CH, Hew, Hw}, an upper bound for the growth of $N(\lam)$ was derived by using a trace estimate in the class of Hilbert-Schmidt operators, obtained with the aid of some inequalities for the Green function of an elliptic operator. In \cite{BS} the asymptotic behavior of eigenvalues was refined by using the Trotter product formula in order to obtain another trace estimate by generalizing the Golden-Thompson inequality, and in \cite{BL} were obtained more terms in the asymptotic expansion by exploiting certain self-similarity of the horns. In Section \S 6 we apply our previous results to this kind of problems in $\R^2$. The main novelty here is the precise order of growth of $N(\lam, \Omega)$ for horns which are not decaying as powers, although is less precise for this kind of horns since the precise constant in the main term is known.

\subsection*{Organization of the paper}
The paper is organized as follows. In Section \S 2 we review some preliminaries results and we introduce the notation that will be used in the rest of the paper. In Section \S 3 we recall the notion of generalized Minkowski dimension and content. In Section \S 4 we estimate the number of eigenvalues of problem \eqref{autovalores} for domains of finite measure. In Section \S 5 we extend this result to infinite measure domains. Finally, in Section \S 6 we estimate the number of eigenvalues of the Laplace operator in two-dimensional horns.

\section{Notations and Preliminaries}

\subsection{Notation}
Throughout the paper, the following notation will be used.

We write $\phi(x)=O(\psi(x))$ when $x\to x_0$ to mean that $|\phi(x)|\leq C\psi(x)$ when $x\to x_0$ for some positive constant $C$. We write $\phi(x)=o(\psi(x))$ when $x \to x_0$ to mean that $\phi(x)/\psi(x) \to 0$ when $x\to x_0$.

Also, we write $\phi(x)\sim \psi(x)$ when $x \to x_0$ to mean that $\phi(x)/\psi(x) \to 1$ when $x \to x_0$, and $\phi(x)\asymp \psi(x)$ when $x \to x_0$ to mean that $c\psi(x)\leq \phi(x)\leq C\psi(x)$ when $x \to x_0$ for some positive constant $c$ and $C$.

\subsection{Eigenvalues of the one dimensional $p-$laplacian}
When $\Omega$ is a interval, in \cite{MANA} the authors obtain a closed formula for the eigenvalues of \eqref{ecu}.

\begin{lema}[See \cite{MANA}, Theorem 3.1] \label{lema1}
Let $ \{ \lambda_{k} \}_{k \in \N} $ be the eigenvalues of \eqref{ecu} in $(0,T)$. Then,
\begin{equation}  \label{cuatro}
\lam_{k}=\frac{\pi_{p}^{p}}{T^{p}}k^{p},
\end{equation}
where $\pi_p$ is given by
$$
\pi_p := 2(p-1)^{1/p}\int_0^1 \frac{ds}{(1-s^p)^{1/p}}.
$$
\end{lema}

From Lemma \ref{lema1} it is easy to see that
$$
N(\lam,(0,T))=\frac{T}{\pi_{p}} \lam^{1/p} + O(1).
$$

The case where $\Omega$ is a disjoint union of intervals, was treated, for instance, in \cite{BOPI}. In that paper, the authors proved the following,
\begin{prop} \label{lema3}
Let $\Omega = \bigcup_{j\in\N} I_j$, where $\{ I_j\}_{j\in\N}$ is a pairwise disjoint family of intervals. Then,
\begin{equation}  \label{cinco}
 N(\lam,\Omega)=\sum_{j=1}^\infty N(\lam, I_j).
\end{equation}
\end{prop}

The following Theorem was proven un \cite{BOPI} and is a suitable generalization of the Dirichlet--Neumann bracketing method of Courant.

\begin{thm}[\cite{BOPI}, Theorem 2.1]\label{Dir-Neu}
Let $U_1, U_2 \in \R^n $ be disjoint open sets such that
$(\overline{U_1 \cup U_2})^{\circ} = U $ and $|U \setminus
U_1 \cup U_2| = 0, $ then
$$
N_D(\lam, U_1 \cup U_2 ) \le N_D(\lam, U) \le N_N(\lam, U) \le N_N(\lam, U_1 \cup U_2).
$$
Here, $N_D(\lam, U)$ (resp., $N_N(\lam, U)$) is the spectral counting function of the Laplace operator in $U$ with Dirichlet boundary conditions on $\partial U$ (resp., with Neumann boundary conditions).
\end{thm}

\subsection{Euler MacLaurin Summation Formula}

We recall the well known summation formula of Euler-MacLaurin, see \cite{KRA} for a proof:

\begin{thm}\label{E-Mc}
Let $f(t)$ be a non negative, continuous and monotonically decreasing function tending to zero when $t \to + \infty$. Then, there exist $C \in \R$, depending only on $f$, such that
\begin{equation}  \label{seis}
\sum_{j=a}^b f(j) = \int_a^b f(t)\, dt + C + O(f(b)),
\end{equation}
when $b \to + \infty$. In particular
\begin{equation}  \label{siete}
\lim_{b \to + \infty}\Big(\sum_{j=a}^b f(j) - \int_a^b f(t)\, dt\Big) = C.
\end{equation}
\end{thm}

\section{Generalized Minkowski content and Minkowski dimension}

\subsection{Minkowski dimension and content}

We denote by $|A|$ the Lebesgue measure of the set $A\subset \R^n$. Let $A_{\varepsilon}$ denote the tubular neighborhood of radius $\varepsilon$ of a set $A\subset \R^n$, i.e.
\begin{equation}  \label{uno}
A_{\varepsilon}=\{ x\in \R^n \colon \mbox{dist}(x,A)\leq \varepsilon \}.
\end{equation}

We recall the classical definition of Minkowski dimension and content (see \cite{Fa, HL, La1, Tr}).

Given $d>0$, the $d-$dimensional upper Minkowski content of $\partial\Omega$ is defined as
\begin{equation}  \label{tres}
M^*(d; \partial \Omega) :=
\limsup_{\varepsilon \to 0^+} \varepsilon^{-(n-d)} |(\partial \Omega)_{\varepsilon} \cap \Omega|.
\end{equation}
Similarly, the $d-$dimensional lower Minkowski content, $M_*(d,\partial\Omega)$, is defined changing the upper by the lower limit in \eqref{tres}.

The Minkowski dimension of $\partial \Omega$ is then defined by
\begin{equation}  \label{dos}
\dim(\partial \Omega) := \inf \{ d \geq 0 \colon M^*(d; \partial \Omega) < \infty\} = \sup \{ d \geq 0\colon M^*(d; \partial \Omega) = \infty\}.
\end{equation}

We will further say that $\partial \Omega$  is $d-$Minkowski measurable if
$$
0 < M_*(d; \partial \Omega) = M^*(d; \partial \Omega) < \infty \qquad \mbox{for some } d>0,
$$
and we call this value $M(d; \partial \Omega)$ the $d-$dimensional Minkowski content of $\partial \Omega$.

\subsection{Dimension functions}
In this paper we will be interested in a suitable generalization of the previous concepts. To this end, given $0<d<1$ we define $G_d$ to be the class of functions $h\colon (0,\infty) \to (0,\infty)$ continuous such that
\begin{itemize}
\item[(H1)] $h$ is stricly increasing and
$$
\lim_{x \to 0^+} h(x)=0,\quad \lim_{x \to \infty} h(x)= \infty.
$$

\item[(H2)] For any $t>0$,
$$
\lim_{x \to 0^+} \frac{h(tx)}{h(x)}=t^d,
$$
uniformly in $t$ on compact subsets of $(0,\infty)$.

\item[(H3)] $h$ is sublinear at $0$, i.e.
$$
\lim_{x \to 0^+} \frac{h(x)}{x}=\infty.
$$
\end{itemize}
One can check that the functions
\begin{equation}  \label{ocho}
h(x)=\frac{x^d}{(\log(\frac{1}{x}+1))^a} \mbox{ and } h(x)=\frac{x^d}{(\log(\log(\frac{1}{x}+1)))^a}
\end{equation}
are in $G_d$ for all $d \in (0,1)$ and $a\geq 0$.

\medskip

\begin{rem} Let $i:(0,\infty)\to(0,\infty)$ be the function $i(x)=x^{-1}$. From now on, given $h \in G_{d}$, we will always let
\begin{equation}  \label{nueve}
g(x):= (h^{-1}\circ i)(x) = h^{-1}(1/x), \quad f(x):= (i\circ h \circ i)(x) = \frac{1}{h(1/x)}.
\end{equation}
\end{rem}

\medskip

With this notations let us now define the generalized Minkowski content and dimension that was introduced by He and Lapidus in \cite{HL}.

\medskip

\begin{defn}
Let $\Omega\subset\R^n$ be an open set with finite Lebesgue measure. Let $h\in G_d$ be a dimension function. The upper $h-$Minkowski content of $\partial\Omega$ is defined by
\begin{equation}\label{content-h}
M^*(h; \partial\Omega) := \limsup_{\ve\to 0^+} \ve^{-n} h(\ve) |(\partial\Omega)_{\ve}\cap \Omega|.
\end{equation}
We define the lower $h-$Minkowski content $M_*(h;\partial\Omega)$ by taking the lower limit in \eqref{content-h}. We further say that $\partial\Omega$ is $h-$Minkowski measurable if
$$
0 < M_*(h;\partial\Omega) = M^*(h;\partial\Omega) < \infty
$$
and denote this value as $M(h;\partial\Omega)$ the $h-$Minkowski content of $\partial\Omega$.
\end{defn}

Let $\Omega$ be an open set in $\R$. Then, $\Omega = \bigcup_{j=1}^{\infty} I_j$, where $I_j$ is an interval of length $l_j$. We can assume that
$$
l_1 \geq l_2 \geq \cdots \geq l_j \geq \cdots > 0
$$

In \cite{HL}, the authors obtained the following relation between the lengths $l_j$ and the Minkowski measurability of $\partial \Omega$:
\begin{thm}
Let $\Omega = \bigcup_{j=1}^{\infty} I_j$. Then, $\partial \Omega$ is $h-$Minkowski measurable if and only if $l_j \sim Lg(j)$. Moreover, in this case, the $h-$Minkowski content of $\partial \Omega$ is given by
$$
M(h;\partial\Omega) = \frac{2^{1-d}}{1-d} L^d.
$$
\end{thm}

Note that $d$ being positive and less than one implies the integrability at infinity of the function $g$, which in turn implies that the Lebesgue measure of the set $\Omega$ is finite. Therefore, the $h-$Minkowski content and dimension are well-defined concepts.

The following proposition, that can be found in \cite{LA}, is a usefull estimate in our arguments in order to compute the constants appearing from the Euler-McLaurin formula.
\begin{prop} \label{proplapidus}
Suppose $h \in G_d$ for some $d \in (0,1)$. Then,
$$
\lim_{x \to \infty} \frac{\int_x^\infty g(u)\, du}{xg(x)} = \frac{d}{1-d}
$$
\end{prop}

\subsection{Nonintegrable Dimension Functions}
We now consider the analogous of the dimension functions defined in the previous subsection to the case $d>1$.

To this end we define the class $G_d$ to be the class of functions $h\colon (0,\infty) \to (0,\infty)$ continuous such that (H1) and (H2) are satisfied and, instead of (H3) we require superlinearity at $0$, i.e.
\begin{equation*}\tag{H3'}
\lim_{x \to 0} \frac{h(x)}{x} = 0.
\end{equation*}

\begin{rem}\label{tilde.h}
As in the previous subsection, we let $i\colon (0,\infty)\to (0,\infty)$ given by $i(x)=x^{-1}$ and
$$
g(x):= (h^{-1}\circ i)(x) = h^{-1}(1/x), \quad f(x):= (i\circ h \circ i)(x) = \frac{1}{h(1/x)}.
$$
\end{rem}

Now we prove an analogous of Proposition \ref{proplapidus} to this case.
\begin{prop} \label{obsinf}
Suppose $h \in G_{d}$ for some $d>1$. Then,
$$
\lim_{x \to \infty} \frac{\int_1^x g(u)\,  du} {xg(x)} = \frac{d}{d-1}
$$
\end{prop}

\begin{proof}
First, we need to show that hypotheses (H1), (H2) and (H3') imply
\begin{equation}\label{homogeneidad.g}
\lim_{x\to\infty} \frac{g(sx)}{g(x)} = s^{-1/d}
\end{equation}
uniformly on $[s_0,\infty)$ for any $s_0>0$ and
\begin{equation}\label{g.infinito}
\lim_{x\to\infty} xg(x) = \infty.
\end{equation}

Equation \eqref{g.infinito} is immediate from (H3'). Now, to prove \eqref{homogeneidad.g} we first observe that it is equivalent to
\begin{equation}\label{h-1}
\lim_{x\to0+} \frac{h^{-1}(sx)}{h^{-1}(x)} = s^{1/d},
\end{equation}
on compact sets of $(0,\infty)$. In order to prove \eqref{h-1}, let us note that (H2) implies
$$
h(sx) = s^d h(x) + o(1),
$$
uniformly on $x$ and in $s\in [0,s_0]$. Then, by the monotonicity of $h$,
$$
h^{-1}(s^d h(x) - \ve) \le sx \le h^{-1}(s^d h(x) + \ve).
$$
Finally, if we call $y=h(x)$ and $t=s^d$,
$$
h^{-1}(ty - \ve) \le t^{1/d} h^{-1}(y) \le h^{-1}(ty + \ve),
$$
which trivially implies \eqref{h-1} and hence \eqref{homogeneidad.g}.

Whith these observations, now the proof of the Proposition follows easily. In fact, by \eqref{g.infinito}, it is enough to prove
$$
\lim_{x \to \infty} \frac{\int_{x_0}^x g(u)\,  du} {xg(x)} = \frac{d}{d-1}
$$
for $x_0$ large enough. Now, by \eqref{homogeneidad.g},
$$
\frac{\int_{x_0}^x g(u)\, du}{xg(x)} = \int_{x_0/x}^1 \frac{g(xs)}{g(x)}\, ds = \int_{x_0/x}^1 s^{-1/d} + o(1)\, ds = \frac{d}{d-1} + o(1).
$$
This fact completes the proof.
\end{proof}

\begin{rem}
Let $\Omega = \bigcup_{j=1}^\infty I_j$ where $I_j$ are disjoint open intervals of lenght $l_j \asymp g(j)$ where $g$ is associated to a function $h\in G_d$ with $d>1$.

In this case, since $g$ is not integrable at infinity, one can check that $|(\partial\Omega)_\ve\cap\Omega| = \infty$ for every $\ve>0$. So, we cannot define the corresponding $h-$Minkowski content or dimension in this case.

Nevertheless, in the computation of the asymptotic behavior of the eigenvalues, we obtain an order of growth for $N(\lambda)$ which depends on $f=(i\circ h\circ i)$.

So, in some sense, $h$ can be considered as certain {\em spectral dimension} for $\partial\Omega$. That is why we refer to $h$ as a {\em nonintegrable dimension function} even though there is no concept of dimension associated to it. See Remark \ref{spectral.dimension} at the end of Section \S 5.
\end{rem}

\section{The finite measure case: $0<d<1$}
An estimate of the number of eigenvalues of the $p-$laplacian equation \eqref{ecu} relies on Lemma \ref{lemaprincipal1} below. This Lemma has been proved in \cite{HL} but we provide here a different proof that will allow us, in the next section, to deal with the infinite measure case.
\begin{lema}\label{lemaprincipal1}
Let $\{l_{j}\}_{j\in\N}$ be an arbitrary nonincreasing positive sequence such that for some $h \in G_{d}$ we have that $l_{j}= g(j)$. Then
$$
\sum_{j=1}^{\infty} [ l_{j} x ] = \sum_{j=1}^{\infty} l_{j}x + \zeta (d)f(x) + o(f(x)), \ \ as \ \ j \to \infty
$$
\end{lema}

\begin{proof}
First, we need to control the difference between $\sum[l_j x]$ and $\sum l_j x$.

To this end, we firs observe that $[l_j x]=0$ if $l_j x<1$. Therefore, the first sum is finite.

Let $J \in \R$ such that $xg(J)=1$. Therefore,
$$
J=g^{-1}\Big(\frac{1}{x}\Big)=\frac{1}{h(1/x)}=f(x).
$$
As $[g(j)x]=0$ if $j>J$, we get
$$
\sum_{j=1}^{\infty} [l_{j}x] = \sum_{j=1}^{J}[g(j)x] = \sum_{j=1}^J g(j) x + O(J).
$$

Observe that this equation immediately gives
$$
\sum_{j=1}^{\infty} [ l_{j} x ] = \sum_{j=1}^{\infty} l_{j}x + O(f(x)).
$$
The rest of the proof will consists in refining the error term.

To improve the remainder estimate, we use Dirichlet's argument for the number of lattice points below the hyperbola: we count the points below the graph of the function $xg(t)$ and below its inverse $g^{-1}(t/x)$, up to the intersection point of these graphs and deleting the size of the square wich we counted twice.

\begin{figure}\label{dibujito}
$$
\begin{array}{cc}
\includegraphics[width=5.7cm]{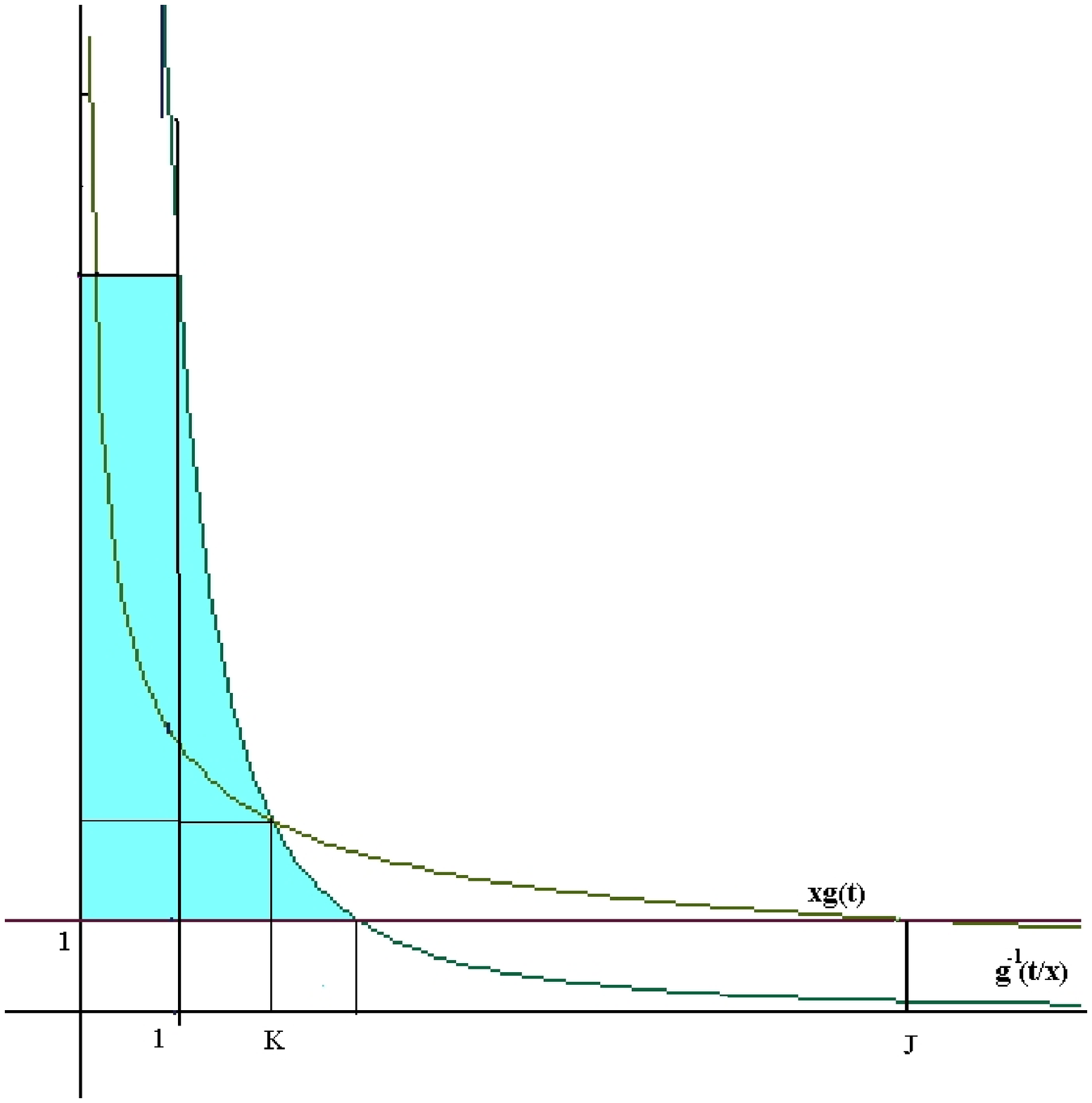} &
\includegraphics[width=5.7cm]{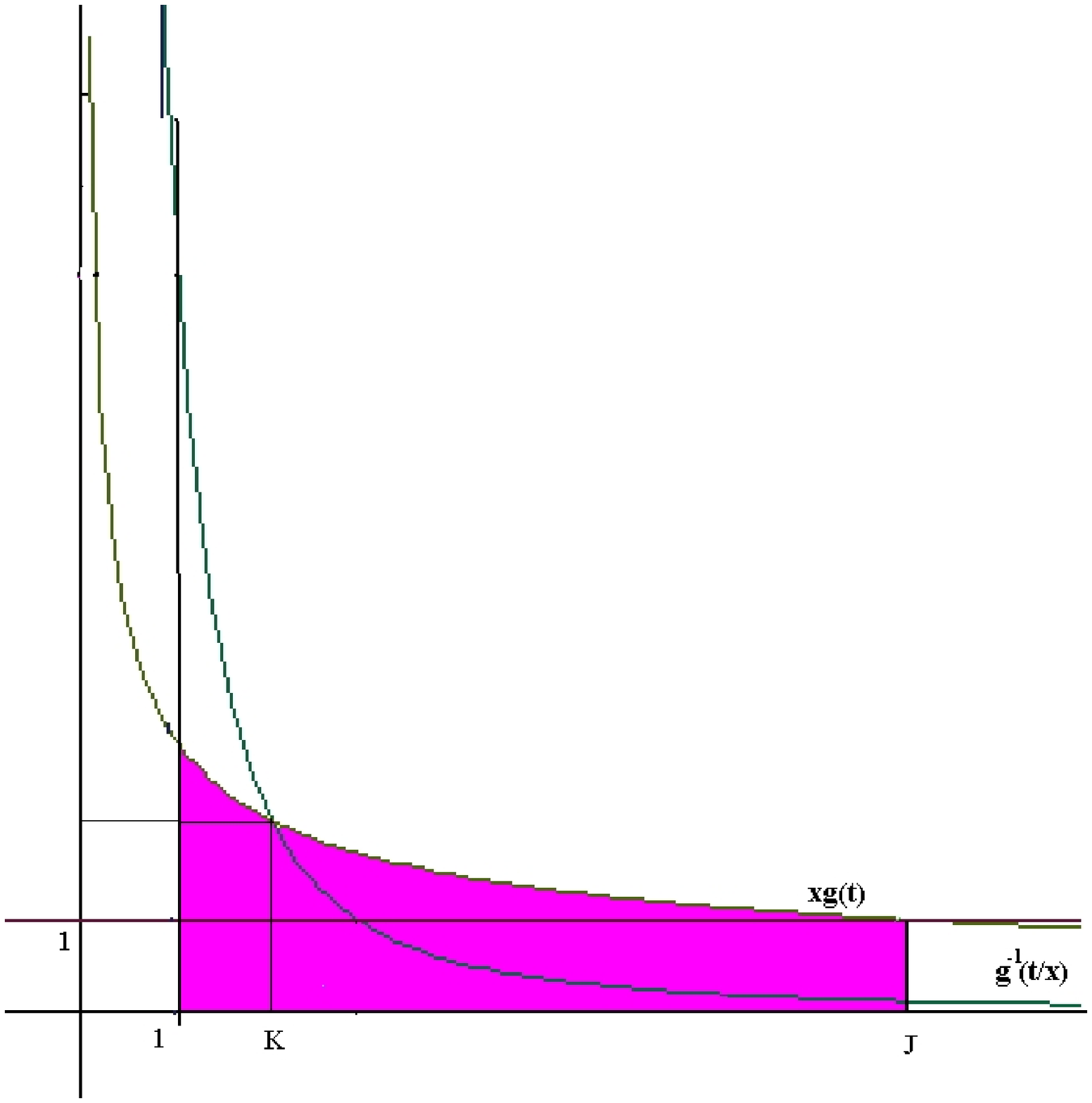}
\end{array}
$$
\caption{Symmetry argument in the proof}
\end{figure}

So, let $K \in \R$ be such that
$$
xg(K) = g^{-1}\Big(\frac{K}{x}\Big) = K
$$
Then
$$
K = g^{-1}\Big(\frac{K}{x}\Big) = \frac{1}{h(\frac{K}{x})} = f\Big(\frac{x}{K}\Big).
$$
By symmetry we have:
\begin{align*}
\sum_{j=1}^{J}[g(j)x] &= \sum_{j=1}^{K}[g(j)x] + \sum_{j=K}^{J}[g(j)x] \\
&= \sum_{j=1}^{K}[g(j)x] + \sum_{j=1}^{K}\Big[g^{-1}\Big(\frac{j}{x}\Big)\Big] - [K]^2 \\
&= \sum_{j=1}^{K}g(j)x + \sum_{j=1}^{K} g^{-1}\Big(\frac{j}{x}\Big) - K^2 + O(K).
\end{align*}
Applying the Euler-McLaurin summation formula \eqref{seis}
\begin{align*}
\sum_{j=1}^{J}[g(j)x] =& \int_{1}^{K}g(t)xdt + A(x) + O(g(K)x) \\
&+ \int_{1}^{K} g^{-1}\Big(\frac{t}{x}\Big)dt + B(x) + O\Big(g^{-1}\Big(\frac{K}{x}\Big)\Big) - K^2 + O(K).
\end{align*}
Clearly
$$
O(K)=O(xg(K))=O\Big(g^{-1}\Big(\frac{K}{x}\Big)\Big).
$$

By symmetry (see Figure 1)
$$
\int_1^J g = \int_1^K g + \int_K^J g = \int_1^K g + \int_1^K g^{-1} - K^2 + J,
$$
then replacing $\int_1^K g(t)x\, dt + \int_1^K g^{-1}\Big(\frac{t}{x}\Big)\, dt$ in the previous equation we have
\begin{equation}\label{f2}
\sum_{j=1}^{J}[g(j)x] = \int_1^J xg(t)\, dt - J + A(x) + B(x) + O\Big(f\Big(\frac{x}{K}\Big)\Big).
\end{equation}

Being the integral convergent, we may write the equation \eqref{f2} as
$$
\sum_{j=1}^{J}[g(j)x] = \int_1^\infty xg(t)\, dt - \int_J^\infty xg(t)\, dt - J + A(x) + B(x) + O\Big(f\Big(\frac{x}{K}\Big)\Big)
$$
and again, by using the Euler-MacLaurin summation formula \eqref{seis}, we obtain
\begin{equation}\label{f3}
\sum_{j=1}^{J}[g(j)x] = \sum_{j=1}^\infty xg(j) - \int_J^\infty xg(t)\, dt - J + B(x) + O\Big(f\Big(\frac{x}{K}\Big)\Big).
\end{equation}

Using that as $x\to\infty$, $K\to\infty$ and by (H2), we obtain $f(x/K) = K \sim f(x)^{1/(1+d)}$. Then
\begin{equation}\label{f4}
\sum_{j=1}^{J}[g(j)x] = x\left(\sum_{j=1}^\infty g(j)\right) - x\int_J^\infty g(t)\, dt - J + B(x) + O(f^{1/(1+d)}(x)).
\end{equation}

To compute the integral we use the Proposition \ref{proplapidus} to obtain
$$
\int_{f(x)}^\infty g(u)\, du = f(x)g(f(x))\Big(\frac{d}{1-d} + o(1)\Big).
$$
Hence, using that $J=f(x)$ and that $g(f(x)) = \frac{1}{x}$ we arrive at
\begin{equation}\label{f5}
x\int_{J}^\infty g(t)\, dt = f(x)\Big(\frac{d}{1-d} +  o(1)\Big),\quad \mbox{as } x \to \infty.
\end{equation}

Replacing in equation \eqref{f4} we obtain
\begin{equation}\label{f6}
\sum_{j=1}^{J}[g(j)x] = x\left(\sum_{j=1}^\infty g(j)\right) - \frac{1}{1-d}f(x) + B(x) + o(f(x)).
\end{equation}

Our last task is to determinate the value of $B(x)$. For $b>1$ fixed, we have,
$$
\sum_{j=1}^b g^{-1}\Big(\frac{j}{x}\Big) - \int_1^b g^{-1}\Big(\frac{t}{x}\Big)\, dt = B(x) + O\Big(g^{-1}\Big(\frac{b}{x}\Big)\Big)
$$
Taking $x$ big enough and remembering that $g^{-1}(t/x)=1/h(t/x)$, $f(x)=1/h(1/x)$,
$$
\sum_{j=1}^b \frac{h(\frac{1}{x})}{h(\frac{j}{x}) h(\frac{1}{x})} - \int_1^b \frac{h(\frac{1}{x})} {h(\frac{t}{x}) h(\frac{1}{x})}\, dt = B(x) + O\Big(g^{-1}\Big(\frac{b}{x}\Big)\Big).
$$
By (H2), for $x$ large we have
$$
\frac{h\Big(\frac{1}{x}\Big)}{h\Big(\frac{t}{x}\Big)} = t^{-d} + o(1).
$$
When $b \to \infty$, as $g^{-1}$ is decreasing, $O(g^{-1}(\frac{b}{x})) \to 0$. Hence,
\begin{equation}\label{f7}
B(x) = f(x)(1+o(1)) \lim_{b \to \infty}\Big( \sum_{j=1}^b j^{-d} - \int_1^b t^{-d} dt \Big),
\end{equation}
or equivalently, $B(x) \sim Cf(x)$ as $x \to +\infty$.
In order to find the constant $C$, we use the next expression for the Riemann zeta function, see \cite{LA}:
$$
\lim_{b \to \infty}\Big( \sum_{j=1}^b j^{-d} - \int_1^b t^{-d} dt \Big) = \zeta(d) - \frac{1}{d-1}.
$$
Hence, replacing in \eqref{f6} the expression $B(x)=f(x)(\zeta(d) - \frac{1}{d-1} + o(1))$ we have
\begin{align*}
\sum_{j=1}^\infty [g(j)x] &= \sum_{j=1}^{J}[g(j)x]\\ &= x\left(\sum_{j=1}^\infty g(j)\right) - \frac{1}{1-d}f(x) + B(x) + o(f(x))\\
&= x\left(\sum_{j=1}^\infty g(j)\right) + \zeta(d)f(x) + o(f(x))
\end{align*}
and the proof is complete.
\end{proof}

Now, we can prove our first theorem:
\begin{thm}\label{caso.finito}
Let $\Omega = \bigcup_{j\in\N} I_j \subset \R$ where $I_j$ are disjoint open intervals. Assume that there exist $d\in (0,1)$ and $h\in G_d$ such that $|I_j|= g(j)$. Then,
$$
N(\lam,\Omega) = \frac{|\Omega|}{\pi_p}\lam^{1/p} + \frac{\zeta(d)}{\pi_p^d}f(\lam^{1/p}) + o(f(\lam^{1/p})) \quad \mbox{as} \quad \lam \to \infty.
$$
\end{thm}

\begin{proof}
As $\Omega = \bigcup_{j\in\N} I_j $ with $|I_j|=g(j)$, from Proposition \ref{lema3},
$$
N(\lam,\Omega) = \sum_{j=1}^\infty \Big[\frac{g(j)}{\pi_p} \lam^{1/p} \Big].
$$
Now the proof follows by a direct application of Lemma \ref{lemaprincipal1} with $x=\lam^{1/p}/\pi_p$. In fact,
$$
N(\lam,\Omega) = \sum_{j=1}^\infty \Big[\frac{g(j)}{\pi_p} \lam^{1/p} \Big] = \frac{|\Omega|}{\pi_p} \lam^{1/p} + \zeta(d) f\Big(\frac{\lam^{1/p}}{\pi_p}\Big) + o(f(\lam^{1/p})),
$$
as we wanted to prove.
\end{proof}

\begin{rem}
Observe that the assumptions of Theorem \ref{caso.finito} implies the length of the intervals $I_j$ to be strictly decreasing. This is not desirable for many applications (for instance, complements of Cantor-type sets).

However, a simple inspection of the arguments show that it suffices to assume that $|I_j|\sim g(j)$. Therefore, for example, complements of Cantor-type sets are included in our result. See \cite{PINA} for the details and also the next section.
\end{rem}

\section{The infinite measure case: $d>1$}

We begin with a couple of lemmas in the spirit of Lemma \ref{lemaprincipal1}.

\begin{lema}\label{lemaprincipal2}
Given $\{l_j\}_{j\in\N}$ a sequence of positive numbers and $h \in G_d$ for some $d>1$. Then, if $l_j \asymp g(j)$, we have
$$
\sum_{j=1}^\infty [l_j x] \asymp f(x) \quad \mbox{as } x \to + \infty.
$$
\end{lema}

\begin{proof}
Since $l_j \asymp g(j)$, there exist two positive constants $c_1,c_2$ such that $c_1 g(j) \leq l_j \leq c_2 g(j)$. Then
$$
c_1 x g(j) - 1 \leq [c_1 x g(j)] \leq [l_j x] \leq [c_2 x g(j)] \leq c_2 x g(j).
$$
So, if we denote $J_i=f(c_ix)$, $i=1,2$ we have that $l_j x < 1$ for $j>J_2$. Then
\begin{equation} \label{form}
\sum_{j=1}^{J_1} c_1 x g(j) - J_1 \leq \sum_{j=1}^{J_2} [l_j x]  = \sum_{j=1}^\infty [l_j x] \leq \sum_{j=1}^{J_2} c_2 x g(j).
\end{equation}

From the summation formula \eqref{E-Mc}, we get
\begin{equation} \label{tetra}
\sum_{j=1}^{J_i} c_i x g(j) = c_i x \int_1^{J_i} g(t)\, dt + Cx + O(xg(J_i))
\end{equation}
Applying Proposition \ref{obsinf}, since $J_i \to \infty$ as, $x \to \infty$
$$
\frac{\int_1^{J_i} g(t)\, dt}{J_i g(J_i)}=\frac{d}{d-1} + o(1).
$$
Also, as $J_i=f(c_i x)$, we have that $c_i x g(J_i)=1$. Moreover, by (H3'), $x = o(f(x))$. Collecting all these facts, we arrive at
$$
\sum_{j=1}^{J_i} c_i x g(j) = \frac{d}{d-1} J_i + o(J_i).
$$
Replacing in \eqref{form} we get
$$
\frac{1}{d-1} J_1 + o(J_1) \le \sum_{j=1}^\infty [l_j x] \le \frac{d}{d-1} J_2 + o(J_2).
$$
Finally, it is easy to see (from (H2)) that $J_i = f(c_i x) \asymp f(x)$ so (1) follows.
\end{proof}

\begin{lema}\label{lemaprincipal3}
Given $\{l_j\}_{j\in\N}$ a sequence of positive numbers and $h \in G_d$ for some $d>1$. Then, if $l_j \sim g(j)$, we have
$$
\sum_{j=1}^\infty [l_j x] = \zeta(d) f(x) + o(f(x)) \quad \mbox{as } x \to + \infty.
$$
\end{lema}

\begin{proof}
Since $l_j \sim g(j)$, for a fixed $\epsilon>0$ there exists $j_0$ such that, for $j>j_0$,
\begin{equation} \label{epsil}
1-\epsilon < \frac{l_j}{g(j)} < 1+ \epsilon.
\end{equation}
From Lemma \ref{lema1} and Proposition \ref{lema3}
\begin{equation} \label{g3}
\sum_{j=1}^\infty [l_j x] = \sum_{j=1}^{j_0} [g(j) x]  + \sum_{j=1}^{j_0} \Big([l_j x] - [g(j) x]\Big) + \sum_{j=j_0+1}^\infty [l_j x].
\end{equation}
Now, from \eqref{epsil} and \eqref{g3} we get
$$
\sum_{j=1}^\infty [(1-\ve)g(j) x] \le \sum_{j=1}^\infty [l_j x] - \sum_{j=1}^{j_0} \Big([l_j x] - [g(j) x]\Big)\le \sum_{j=1}^\infty [(1+\ve)g(j) x]
$$
Now, if $K_\pm$ is such that
$$
(1\pm\ve)g(K_\pm)x = g^{-1}(K_\pm/x(1\pm\ve)) = K_\pm,
$$
arguying as in Lemma \ref{lemaprincipal1}, we arrive at
$$
\sum_{j=1}^\infty [(1\pm \ve)g(j)x] = \sum_{j=1}^{K_\pm} (1\pm\ve) g(j) x + \sum_{j=1}^{K_\pm} g^{-1}(j/x(1\pm\ve)) - K_\pm^2 + O(K_\pm).
$$

Applying the Euler-McLaurin summation formula \eqref{E-Mc}, we get
\begin{align*}
\sum_{j=1}^\infty [(1\pm \ve)g(j)x] =& \int_1^{K_\pm} (1\pm\ve) g(t)x\, dt + \int_1^{K_\pm} g^{-1}(t/x(1\pm\ve))\, dt\\
& + A(x) + B(x) - K_\pm^2 + O(K_\pm),
\end{align*}
where $A(x) = C(1\pm\ve)x$ and $B(x)$ are the constants from the Euler-McLaurin formula \eqref{E-Mc} for $(1\pm\ve) g(t)x$ and $g^{-1}(t/x(1\pm\ve))$ respectively.

Again, as in Lemma \ref{lemaprincipal1}
$$
\int_1^{K_\pm} g + \int_1^{K_\pm} g^{-1} = \int_1^{J_\pm} g + K_\pm(K_\pm -1) - J_\pm,
$$
where $J_\pm$ is given by $(1\pm\ve)xg(J_\pm) = 1$.

Therefore, we arrive at
$$
\sum_{j=1}^\infty [(1\pm \ve)g(j)x] = \int_1^{J_\pm} (1\pm\ve)xg(t)\, dt + A(x) + B(x) - J_\pm + O(K_\pm).
$$
Applying now Proposition \ref{obsinf} and the definition of $J_\pm$ we obtain
\begin{align*}
\sum_{j=1}^\infty [(1\pm \ve)g(j)x] =& (1\pm\ve)x J_\pm g(J_\pm) \Big(\frac{d}{d-1} + o(1)\Big)\\
& + A(x) + B(x) - J_\pm + O(K_\pm)\\
=& J_\pm \Big(\frac{1}{d-1} + o(1)\Big) + A(x) + B(x) + O(K_\pm)\\
=& \frac{1}{d-1} f((1\pm\ve)x) + B(x) + o(f(x)),
\end{align*}
where we have used that $A(x) = C(1\pm\ve)x$, $x=o(f(x))$ and $K_\pm = f(x(1\pm\ve)/K_\pm) = o(f(x))$.

It remains to estimate $B(x)$, but this follows exactly as in the proof of the finite measure case, Proposition \ref{lemaprincipal1}. So
$$
B(x) = f((1\pm\ve)x)(1+o(1)) \lim_{b \to \infty}\Big(\sum_{j=1}^b j^{-d} - \int_1^b t^{-d}\, dt \Big).
$$
In this case, both terms are convergent, and we easily get
$$
B(x) = \Big(\zeta(d) - \frac{1}{d-1}\Big) f((1\pm\ve)x) + o(f(x)).
$$
Hence, we finally get
$$
\sum_{j=1}^\infty [(1\pm \ve)g(j)x] = \zeta(d) f((1\pm\ve)x) + o(f(x)).
$$
As $\ve>0$ is arbitrary, the proof follows.
\end{proof}

Now, we can prove our second theorem:
\begin{thm}
Let $\Omega = \bigcup_{j\in\N} I_j$, and $h \in G_{d}$ for some $d > 1$. Then
\begin{enumerate}
\item  if $|I_j|_1 \asymp g(j)$, we have
$$
N(\lam,\Omega)=O(f(\lam^{1/p})) \quad \mbox{as } \lam \to + \infty.
$$

\item  if $|I_j|_1 \sim g(j)$, we have
$$
N(\lam,\Omega) = \frac{\zeta(d)}{\pi_p^d} f(\lam^{1/p}) + o(f(\lam^{1/p})) \quad \mbox{as } \lam \to + \infty.
$$
\end{enumerate}
\end{thm}

\begin{proof}
The proofs follow from Lemmas \ref{lemaprincipal2} and \ref{lemaprincipal3} replacing $x$ by $\lam^{1/p}/\pi_p^{1/p}$.
\end{proof}

We close this section with the following estimate for the eigenvalues.

\begin{cor}Let $h\in G_d$ for some $d>1$ and let $\Omega = \bigcup_{j\in\N} I_j$ be such that $|I_j|\sim g(j)$.
Let $\{\lam_k\}_{k\in \N}$ be the sequence of eigenvalues of problem \eqref{ecu} in $\Omega$. Then,
$$\lam_k \sim \Big[ g\Big(\frac{\pi_p^d k}{\zeta(d)}\Big)\Big]^{-p}.$$
\end{cor}

\begin{proof} Since
$$
k = N(\lam_k,\Omega) \sim \frac{\zeta(d)}{\pi_p^d} f(\lam_k^{1/p}) = \frac{\zeta(d)}{\pi_p^d} g^{-1}(\lam_k^{-1/p}),
$$
we get
$$
\Big[ g\Big(\frac{\pi_p^d k}{\zeta(d)}\Big)\Big]^{-p}   \sim  \lam_k
$$
and the proof is finished.
\end{proof}

\begin{rem}\label{spectral.dimension}
Let us note that, for $h(t)=t^d$, we have that $g(t) = t^{-1/d}$, so
$$
\lam_k  \sim \Big(\frac{\pi_p^d k}{\zeta(d)}\Big)^{p/d}= \frac{\pi_p^p k^{p/d}}{\zeta(d)^{p/d}}.
$$

For $p=2$, the eigenvalues of the Laplace operator with Dirichlet boundary condition in any bounded open set $U \subset \R^N$ satisfy
$$
\lam_k \sim c k^{2/N}.
$$
Hence, seems natural to consider $h$ as a {\em spectral dimension} for $\partial\Omega$ despite the fact that
$\Omega = \bigcup_{j\in\N} I_j \subset \R$ and $d>1$.
\end{rem}

\section{Two--dimensional horns }

For simplicity, we only consider here two dimensional domains.  First we derive a simple proof of the upper bound for the eigenvalue counting function of the Laplace operator on horns. Then, we give a lower bound with the same order of growth although with a different constant in the leading term.

Let $h\in G_d$, with $d>1$, and $g(x) = h^{-1}(1/x)$.  Let $\Omega\subset \R^2$ be defined as
$$
\Omega = \{ (x,y) \in \R^2 :  x\ge 1; \; |y| \le g(x) \}.
$$
Clearly, the measure of $\Omega$ is infinite.

Let us consider the eigenvalue problem
\begin{equation}\label{quasib}
\begin{cases}
-\Delta u = \lam u & \mbox{in } \Omega \\
u = 0 &  \mbox{on } \partial \Omega.
\end{cases}
\end{equation}
Since $g(x)\searrow 0$ as $x\nearrow \infty$, the domain is quasibounded, namely,
$$
\lim_{|x|\to \infty} d(x, \R^2\setminus \Omega) = 0,
$$
and the spectrum is discrete, consisting of a sequence of eigenvalues $0<\lam_1 < \lam_2 \le \cdots \nearrow \infty$, repeated according their multiplicity.

We want to estimate the order of growth of
$$
N(\lam, \Omega) = \#\{n\in\N\colon \lam_n \le \lam\}.
$$
To this end, let us introduce a family of rectangles $\{Q^j\}_{j\in\N}$, and an open set $V$ such that $\Omega \subset V$:
$$
Q^j = [j,j+1]\times[-g(j), g(j)], \qquad V =  \Big(\bigcup_{j=1}^\infty Q^j\Big)^{\circ}.
$$
Also, the set $V$ is quasibounded and the spectrum of the Laplace operator in $V$ is a sequence $\mu_1 \le \mu_2 \le \cdots \nearrow \infty$, repeated according their multiplicity. Moreover, the monotonicity of eigenvalues respect to the domain gives
$$
\mu_n \le \lam_n, \qquad n\ge 1.
$$

We have the following inclusions of Sobolev spaces:
$$
H_0^1(\Omega) \subset  H_0^1(V) \subset \bigoplus_{j=1}^{\infty}  H_*^1(Q^j),
$$
where
$$
H_*^1(Q^j) = \{ u \in H^1(Q^j) : u(x, \pm g(j)) = 0\}.
$$

We can compute by separation of variables the eigenfunctions and eigenvalues of the Laplace operator in each $Q^j$ with mixed boundary conditions. We get
$$
\lam_{h,k}^{Q^j} = h^2\pi^2 + \frac{k^2 \pi^2}{4g(j)^2}, \quad u_{h,k}^{Q_j}(x,y)=\cos(h\pi y)\sin(k\pi y/2g(j)), \quad h\ge 0, k\ge 1.
$$
Hence, we define the eigenvalue counting function
$$
N_{mixed}(\lam, Q^j)= \#\Big\{(h,k) : h^2\pi^2 + \frac{k^2 \pi^2}{4g(j)^2} \le \lam, \quad h\ge 0, k\ge 1\Big\}.
$$
Let us note that we can assign to each eigenvalue a lattice point $(h,k)$ with $h>0$ and the square $(h-1, h]\times (k-1,k]$, and the number of eigenvalues with $h=0$ is $[2g(j)\lam^{1/2}/\pi]$. By using the area of the ellipse which contains those squares, we get
\begin{equation}\label{mixed}
N_{mixed}(\lam, Q^j) \le \frac{g(j)}{2\pi^2}\lam + \frac{2g(j)}{\pi}\lam^{1/2} = g(j) \Big(\frac{\lam}{2\pi^2} + \frac{2\lam^{1/2}}{\pi}\Big).
\end{equation}

Now, the Dirichlet-Neumann bracketing \eqref{Dir-Neu} together with Proposition \ref{lema3} implies
$$
N(\lam, \Omega)\le \sum_{j=1}^\infty N_{mixed}(\lam, Q^j),
$$
but we cannot replace the previous bound yet.
Let us note that $N_{mixed}(\lam, Q^j) = 0$ if
$$
\lam_{0,1}^{Q^j} = \frac{ \pi^2}{4g(j)^2} > \lam,
$$
i.e., for $j > g^{-1}( \pi/2\lam^{1/2}) = f(2\lam^{1/2}/\pi)$. Hence, by using the estimate \eqref{mixed},  the Euler-McLaurin formula \eqref{E-Mc} and Proposition  \ref{obsinf}, we obtain
\begin{align*}
N(\lam, \Omega) \le &  \sum_{j=1}^{f(2\lam^{1/2}/\pi)}N_{mixed}(\lam, Q^j) \\
\le &  \sum_{j=1}^{f(2\lam^{1/2}/\pi)} g(j) \Big(\frac{\lam}{2\pi^2} + \frac{2\lam^{1/2}}{\pi}\Big)\\
= & \Big(\frac{\lam}{2\pi^2} + \frac{2\lam^{1/2}}{\pi}\Big) \left( \int_{1}^{f(2\lam^{1/2}/\pi)} g(t) dt
+ A + O\Big( g\Big(f\Big( \frac{2\lam^{1/2}}{\pi} \Big)\Big)\Big) \right) \\
= & \Big(\frac{\lam}{2\pi^2} + \frac{2\lam^{1/2}}{\pi}\Big) f\Big( \frac{2\lam^{1/2}}{\pi}\Big) g\Big( f\Big( \frac{2\lam^{1/2}}{\pi}\Big)\Big)
\Big( \frac{d}{d-1} + o(1)\Big) + O( \lam^{1/2} ) \\
= & \Big(\frac{\lam}{2\pi^2} + \frac{2\lam^{1/2}}{\pi}\Big) f\Big( \frac{2\lam^{1/2}}{\pi}\Big)  \frac{\pi}{2\lam^{1/2}}
\Big( \frac{d}{d-1} + o(1)\Big) + O( \lam^{1/2} ) \\
= &  \frac{d}{4\pi(d-1)}  \lam^{1/2} f\Big( \frac{2\lam^{1/2}}{\pi}\Big)
+ o\left(\lam^{1/2} f\Big( \frac{2\lam^{1/2}}{\pi}\Big)\right).
\end{align*}

Therefore, we have proved the following Theorem:
\begin{thm}\label{upper-quasib}
Let $h\in G_d$, with $d>1$, and $\Omega\subset \R^2$ be defined as
$$
\Omega = \{ (x,y) \in \R^2 :  x\ge 1; \; |y| \le g(x) \}.
$$
Then, the eigenvalue counting function of the eigenvalue problem \eqref{quasib}
satisfies
$$
N(\lam, \Omega) \le \frac{d}{d-1} \lam^{1/2} f\Big( \frac{2\lam^{1/2}}{\pi}\Big) + o\left(\lam^{1/2} f\Big( \frac{2\lam^{1/2}}{\pi}\Big)\right).
$$
\end{thm}

\begin{rem}
When $h(t)= t^d$ with $d>1$, then $g(t)= t^{-1/d}$ and $f(t) = t^{d}$. So, we have
$$
N(\lam, \Omega) \le  \frac{d}{d-1} \Big(\frac{2}{\pi}\Big)^{d} \lam^{\frac{d+1}{2}}
+ o(\lam^{\frac{d+1}{2}}).
$$
Following \cite{BL}, the order of growth cannot be improved, since this is the right one for horn-shaped domains.
\end{rem}

In much the same way we prove the following lower bound:

\begin{thm}\label{lower-quasib}
Let $h\in G_d$, with $d>1$, and $\Omega\subset \R^2$ be defined as
$$
\Omega = \{ (x,y) \in \R^2 :  x\ge 1; \; |y| \le g(x) \}.
$$
Then, the eigenvalue counting function of the eigenvalue problem \eqref{quasib} satisfies
$$
N(\lam, \Omega) \ge \frac{1}{d-1} \frac{ \lam^{1/2}}{\pi} f\Big( \frac{\lam^{1/2}}{2\pi}\Big) + o\left(\lam^{1/2} f\Big( \frac{ \lam^{1/2}}{2\pi}\Big)\right).
$$
\end{thm}

\begin{proof}
As before, let us introduce a family of rectangles $\{Q_j\}_{j\in\N}$ and $U\subset \Omega$, where
$$
Q_j = [j,j+1]\times[-g(j+1), g(j+1)], \qquad  U = \Big(\bigcup_{j=1}^\infty Q_j\Big)^{\circ}.
$$ 
Then,
$$
\bigoplus_{j=1}^{\infty} H_0^1(Q_j) \subset H_0^1(U),
$$
and the Dirichlet-Neumann bracketing \eqref{Dir-Neu} together with Proposition \ref{lema3} implies
$$
\sum_{j=1}^\infty N_D(\lam, Q_j) \le N(\lam, \Omega).
$$

The eigenfunctions and eigenvalues of the Laplace operator in $Q_j$ with Dirichlet boundary conditions are
$$
\lam_{h,k}^{Q_j} = h^2\pi^2 + \frac{k^2 \pi^2}{4g(j+1)^2}, \quad u_{h,k}^{Q_j}(x,y)=\sin(k\pi x/2g(j))\sin(h\pi y), \quad h,k \ge 1.
$$
Therefore, the counting function $N_D(\lam, Q_j)$ is
$$
N_D(\lam, Q_j) = \#\Big\{(h,k) : h^2\pi^2 + \frac{k^2 \pi^2}{4g(j+1)^2}\le \lam, \quad h,k \ge 1\Big\}.
$$

Let us assign to each eigenvalue the lattice point $(h,k)$ with $h,k \ge1$, and the square $Q_{h,k} = [h, h+1)\times [k,k+1)$. Hence,
$$
N_D(\lam, Q_j) = \Big| \Big( \bigcup_{\lam_{h,k}^{Q_j}\le \lam} Q_{h,k}\Big)\Big|.
$$

Clearly,
$$
N_D(\lam, Q_j) \ge \frac{g(j)\lam}{2\pi^2} - \frac{\lam^{1/2}}{\pi} - \frac{2g(j)\lam^{1/2}}{\pi} -1,
$$
since in the first quadrant, the  ellipse of semi-axes $\lam^{1/2}/\pi$ and $ 2g(j)\lam^{1/2}/\pi$ is covered by the squares $Q_{h,k}$ and the rectangles  $[0,1)\times[0,\lam^{1/2})$, $[0,[2g(j)\lam^{1/2}]+1)\times[0,1)$.

We consider only $j\le f(\lam^{1/2}/2\pi)$ (if not, $\frac{g(j)\lam}{2\pi^2} - \frac{\lam^{1/2}}{\pi} <0$, and $N_D(\lam, Q_j)$ is nonnegative) and we get
$$
N(\lam, \Omega)\ge \sum_{j=1}^\infty N_D(\lam, Q_j) \ge\sum_{j=1}^{f(\lam^{1/2}/2\pi)}\frac{g(j)\lam}{2\pi^2} - \frac{\lam^{1/2}}{\pi} - \frac{2g(j)\lam^{1/2}}{\pi} -1.
$$

Finally, as in the previous proof,
\begin{align*}
N(\lam, \Omega)\ge & \sum_{j=1}^{f(\lam^{1/2}/2\pi)}\frac{g(j)\lam}{2\pi^2}  - f\Big(\frac{\lam^{1/2}}{2\pi}\Big) \frac{\lam^{1/2}}{\pi} + O\Big( \sum_{j=1}^{f(\lam^{1/2}/2\pi)}2g(j)\lam^{1/2} \Big)  \\
=  & \frac{\lam}{2\pi^2} f\Big( \frac{\lam^{1/2}}{2\pi}\Big) g\Big( f\Big( \frac{\lam^{1/2}}{2\pi}\Big)\Big) \Big( \frac{d}{d-1} + o(1)\Big)  - \frac{\lam^{1/2}}{\pi} f\Big(\frac{\lam^{1/2}}{2\pi}\Big) \\
& + o\Big( \lam^{1/2} f\Big( \frac{\lam^{1/2}}{2\pi}\Big) \Big)  \\
=  & \frac{\lam}{2\pi^2} f\Big( \frac{\lam^{1/2}}{2\pi}\Big) \frac{2\pi}{\lam^{1/2}}
\Big( \frac{d}{d-1} + o(1)\Big)  - \frac{\lam^{1/2}}{\pi} f\Big(\frac{\lam^{1/2}}{2\pi}\Big) + o\Big( \lam^{1/2} f\Big( \frac{\lam^{1/2}}{2\pi}\Big) \Big) \\
=  &  \frac{\lam^{1/2}}{\pi(d-1)} f\Big(\frac{\lam^{1/2}}{2\pi}\Big)
+ o\Big( \lam^{1/2} f\Big( \frac{\lam^{1/2}}{2\pi}\Big) \Big)
\end{align*}
and the proof is finished.
\end{proof}

\begin{rem} From Theorems \ref{upper-quasib} and \ref{lower-quasib} we obtain that
$$
c \lam^{1/2} f\Big(\frac{\lam^{1/2}}{2\pi}\Big) \le N(\lam, \Omega) \le C \lam^{1/2} f\Big(\frac{2\lam^{1/2}}{\pi}\Big),
$$
for horn-shaped domains
$$
\Omega = \{ (x,y) \in \R^2 :  x\ge 1; \; |y| \le g(x) \},
$$
with $f(x) = g^{-1}(1/x)$, and $g$ monotonically decreasing continuous function.

Observe that, as $h$ satisfies (H2), we have
$$
N(\lam, \Omega) \asymp \lam^{1/2} f(\lam^{1/2}).
$$
This result improves the upper bounds obtained in \cite{CH, Hew, Hw}, which only gives an upper bound for $N(\lam, \Omega)$ whenever $g(x) = x^{-1/d}$.

It would be desirable to obtain a better knowledge of the asymptotic behavior, namely, $N(\lam, \Omega) \sim
c \lam^{1/2} f(\lam^{1/2})$ (for certain constant $c$) as in \cite{BS}, and even a second term as in \cite{BL}. However, without imposing more restrictions on the functions $h$ or $g$, we believe that this cannot be possible, since the main term can oscillate, as the following one--dimensional example suggest. This example is borrowed from \cite{PINA}.
\end{rem}

\begin{exam} Let $\Omega = \bigcup_{k\in\N} \Omega_k$, where $\Omega_k$ consist of $m^k$ intervals of lengths $n^{1-k}$, for $m>n$. Then, the spectral counting function of problem \eqref{ecu} satisfies
$$
N(\lam, \Omega)=  \frac{\lam^{d/p}}{m} s(\log(\lam)) -  O(
\lam^{1/p}),
$$
where $s(\log(\lam))$ is a bounded  periodic function,  and
$d=\frac{\log(m)}{\log(n)}$.
\end{exam}

\begin{proof} Since
$$
N(\lam, \Omega) = \sum^{\infty}_{j=1} m^j \left[ \frac{\lam^{1/p}}{\pi_p n^{j-1}} \right] =
\sum^{\infty}_{-\infty}  m^j \left[ \frac{\lam^{1/p}}{n^{j-1}\pi_p} \right] - O(\lam^{1/p}).
$$
By changing variables,
$$  \label{vary} k = \frac{\log(\lam^{1/p})- \log(\pi_p)}{\log(n)},
$$
we get $n^k = \lam^{1/p}/\pi_p$ and $m^k = (\lam^{1/p}/\pi_p)^d$, for $d=\frac{\log(m)}{\log(n)}$, and we obtain
$$
N(\lam, \Omega) = \frac{\lam^{d/p}}{m} \sum^{\infty}_{j=-\infty} m^{j-k} [n^{y-k}] - O(\lam^{1/p}) = \frac{\lam^{d/p}}{m}s(\log(\lam))- O(\lam^{1/p})
$$
and, as  $j-(k+1) = (j+1) - k$,  $s(\log(\lam))$ is a periodic function with period equal to one.
\end{proof}

This example can be extended to $\R^2$, by defining $\Omega = \bigcup_{k\in\N} \Omega_k$, where $\Omega_k$ consists of $m^{k}$ disjoints squares of sides $n^{1-k}$. When $\Omega$ has finite measure, similar examples were considered in \cite{FV, LV, BLe}, where oscillating second term were obtained for the spectral counting function of the Laplace operator in $\Omega$ with Dirichlet boundary
conditions in the boundary of each square. It is not difficult to extend those arguments to the infinite measure case (that is, $m^2>n$), to obtain in this way a quasibounded set with an oscillating main term. However, the set obtained in this way is not a horn.

\section*{Acknowledgements}  Supported by grant X078 from Universidad de Buenos Aires, and grants PICT 06-835, PICT 06-290 from ANPCyT. J. Fern\'andez Bonder and J.P. Pinasco are members of CONICET. A.M. Salort is a fellow of CONICET.


\begin{thebibliography}{10}

\bibitem{CH}
Colin Clark and Denton Hewgill.
\newblock One can hear whether a drum has finite area.
\newblock {\em Proc. Amer. Math. Soc.}, 18:236--237, 1967.

\bibitem{MANA}
Pavel Dr{\'a}bek and Ra{\'u}l Man{\'a}sevich.
\newblock On the closed solution to some nonhomogeneous eigenvalue problems
  with {$p$}-{L}aplacian.
\newblock {\em Differential Integral Equations}, 12(6):773--788, 1999.

\bibitem{Fa}
Kenneth Falconer.
\newblock {\em Fractal geometry}.
\newblock John Wiley \& Sons Ltd., Chichester, 1990.
\newblock Mathematical foundations and applications.

\bibitem{BOPI}
Juli{\'a}n Fern{\'a}ndez~Bonder and Juan~Pablo Pinasco.
\newblock Asymptotic behavior of the eigenvalues of the one-dimensional
  weighted {$p$}-{L}aplace operator.
\newblock {\em Ark. Mat.}, 41(2):267--280, 2003.

\bibitem{FV}
Jacqueline Fleckinger-Pell{\'e} and Dmitri~G. Vassiliev.
\newblock An example of a two-term asymptotics for the ``counting function'' of
  a fractal drum.
\newblock {\em Trans. Amer. Math. Soc.}, 337(1):99--116, 1993.

\bibitem{HL}
Christina~Q. He and Michel~L. Lapidus.
\newblock Generalized {M}inkowski content, spectrum of fractal drums, fractal
  strings and the {R}iemann zeta-function.
\newblock {\em Mem. Amer. Math. Soc.}, 127(608):x+97, 1997.

\bibitem{Duke}
Dennis~A. Hejhal.
\newblock The {S}elberg trace formula and the {R}iemann zeta function.
\newblock {\em Duke Math. J.}, 43(3):441--482, 1976.

\bibitem{Hew}
Denton Hewgill.
\newblock On the eigenvalues of a second order elliptic operator in an
  unbounded domain.
\newblock {\em Pacific J. Math.}, 51:467--476, 1974.

\bibitem{Hw}
Denton~E. Hewgill.
\newblock On the growth of the eigevalues of an elliptic operator in a
  quasibounded domain.
\newblock {\em Arch. Rational Mech. Anal.}, 56:367--371, 1974/75.

\bibitem{Kac}
Mark Kac.
\newblock Can one hear the shape of a drum?
\newblock {\em Amer. Math. Monthly}, 73(4, part II):1--23, 1966.

\bibitem{KRA}
Ekkehard Kr{\"a}tzel.
\newblock {\em Lattice points}, volume~33 of {\em Mathematics and its
  Applications (East European Series)}.
\newblock Kluwer Academic Publishers Group, Dordrecht, 1988.

\bibitem{La1}
Michel~L. Lapidus.
\newblock Fractal drum, inverse spectral problems for elliptic operators and a
  partial resolution of the {W}eyl-{B}erry conjecture.
\newblock {\em Trans. Amer. Math. Soc.}, 325(2):465--529, 1991.

\bibitem{LA}
Michel~L. Lapidus and Carl Pomerance.
\newblock The {R}iemann zeta-function and the one-dimensional {W}eyl-{B}erry
  conjecture for fractal drums.
\newblock {\em Proc. London Math. Soc. (3)}, 66(1):41--69, 1993.

\bibitem{LV}
Michael Levitin and Dmitri Vassiliev.
\newblock Spectral asymptotics, renewal theorem, and the {B}erry conjecture for
  a class of fractals.
\newblock {\em Proc. London Math. Soc. (3)}, 72(1):188--214, 1996.

\bibitem{PINA}
Juan~Pablo Pinasco.
\newblock Asymptotic of eigenvalues and lattice points.
\newblock {\em Acta Math. Sin. (Engl. Ser.)}, 22(6):1645--1650, 2006.

\bibitem{BS}
Barry Simon.
\newblock Nonclassical eigenvalue asymptotics.
\newblock {\em J. Funct. Anal.}, 53(1):84--98, 1983.

\bibitem{Tr}
Claude Tricot, Jr.
\newblock Two definitions of fractional dimension.
\newblock {\em Math. Proc. Cambridge Philos. Soc.}, 91(1):57--74, 1982.

\bibitem{BLe}
M.~van~den Berg and M.~Levitin.
\newblock Functions of {W}eierstrass type and spectral asymptotics for iterated
  sets.
\newblock {\em Quart. J. Math. Oxford Ser. (2)}, 47(188):493--509, 1996.

\bibitem{BL}
M.~van~den Berg and M.~Lianantonakis.
\newblock Asymptotics for the spectrum of the {D}irichlet {L}aplacian on
  horn-shaped regions.
\newblock {\em Indiana Univ. Math. J.}, 50(1):299--333, 2001.

\end{thebibliography}
\end{document}